\documentstyle[12pt,amssymb,amstex]{amsart}

\numberwithin{equation}{section}

\newtheorem{theorem}{Theorem}[section]

\newtheorem{lemma}[theorem]{Lemma}

\def\bn{{\mathbb N}}

\def\<{\langle}
\def\>{\rangle}

\def\1{\mathbf{1}}

\begin{document}

\title[Solution of variational inequality on fixed points set]
{A Solution of variational inequality problem for a finite family of nonexpansive mappings in Hilbert spaces}

\author{Farrukh Mukhamedov}
\address{Farrukh Mukhamedov\\
Department of Computational \& Theoretical Sciences \\
Faculty of Sciences, International Islamic University Malaysia\\
P.O. Box, 141, 25710, Kuantan\\
Pahang, Malaysia} \email{{\tt far75m@@yandex.ru}}

\author{Mansoor Saburov}
\address{Mansoor Saburov\\
Department of Computational \& Theoretical Sciences \\
Faculty of Science, International Islamic University Malaysia\\
P.O. Box, 141, 25710, Kuantan\\
Pahang, Malaysia} \email{{\tt msaburov@@gmail.com}}

\begin{abstract}
In this paper we prove the strong convergence of the explicit
iterative process to a common fixed point of the finite family of nonexpansive
mappings defined on Hilbert space, which solves the the variational inequality on the fixed points set.  \vskip 0.3cm \noindent {\it
Mathematics Subject Classification}: 46B20; 47H09; 47H10\\
{\it Key words}: Nonexpansive mappings, explicit iteration process; common fixed point; variational inequality.
\end{abstract}

\maketitle

\section{Introduction}
An iterative approximation of fixed points and zeros of nonlinear operators has been
studied extensively by many authors to solve nonlinear operator equations as well
as variational inequality problems (see e.g., \cite{[14],[15]}, \cite{[20]}, \cite{[25]}-\cite{[31]}, \cite{[34]}-\cite{[40]}).
A very important class of mappings is nonexpansive mappings. In particulary, iterative approximation
of fixed points of nonexpansive mappings is an important subject in nonlinear functional analysis, and has many applications in various fields of mathematics, physics etc
(see e.g., \cite{[12],[33],[50]}).

Let $H$ be a real Hilbert space with inner product $\langle\cdot,\cdot\rangle,$ and
induced norm $\|\cdot\|.$ Let $T:H\to H$ be a mapping. Denote by $F(T)$ the set of fixed points of $T$,
that is, $F(T) =\{x\in H: Tx = x\}$. Now let us recall some well-known definitions, in which will be used  in this paper.

A mapping $f:H\to H$ is said to be \emph{contraction} if there exists $\alpha\in[0,1)$ such that for all $x,y\in H,$
\[\|f(x)-f(y)\|\le \alpha\|x-y\|.\]
A mapping $T:H\to H$ is called \emph{nonexpansive} if for all $x,y\in H,$
\[\|Tx-Ty\|\le \|x-y\|.\]
A mapping $A:H\to H$ is called \emph{$\eta-$strongly monotone} if for all $x,y\in H,$
\[\Bigl\langle x-y,Ax-Ay\Bigr\rangle \ge \eta\|x-y\|^2.\]
A mapping $A:H\to H$ is called \emph{$k-$Lipschitzian} if for all $x,y\in H,$
\[\|Ax-Ay\|\le k\|x-y\|.\]

Iterative methods for nonexpansive mappings are now also applicable in solving
convex minimization problems (see, for example, \cite{[47]} and references therein). Let
$K$ be a closed convex nonempty subset of $H$ and $T: K\to K$ be a nonexpansive mapping such that
$F(T)\neq\emptyset.$ Given $u\in K$ and a real sequence $\{\alpha_n\}_{n=1}^\infty$ in the interval $(0,1),$ starting
with an arbitrary initial $x_0\in K,$ let a sequence $\{x_n\}_{n=1}^\infty$ be defined by
\begin{equation}\label{xnalphanTxn}
x_{n+1} = \alpha_nu + (1-\alpha_n)Tx_n, \quad n\in\bn.
\end{equation}
Under appropriate conditions on the iterative parameter $\{\alpha_n\},$ it has been shown by
Halpern \cite{[18]}, Lions \cite{[23]}, Wittmann \cite{[43]} and Bauschke \cite{[2]} that $\{x_n\}_{n=1}^\infty$ converges
strongly to $P_{F(T)}u,$ the projection of $u$ to the fixed point set, $F(T)$ of $T.$ This means
that the limit of the sequence $\{x_n\}_{n=1}^\infty$ solves the following minimization
problem:
\begin{center}
find $x^{*}\in F(T)$ such that $\|x^{*}-u\|=\min\limits_{x\in F(T)}\|x-u\|.$
\end{center}
A typical minimization problem is to minimize a quadratic function over the set of
fixed points of nonexpansive mappings in a real Hilbert space and such problems
go thus: find $x^{*}\in F(T)$ such that
\begin{equation}
\frac12\Bigl\langle Ax^{*},x^{*}\Bigr\rangle-\Bigl\langle x^{*},u\Bigr\rangle=\min\limits_{x\in F(T)}\left(\frac12\Bigl\langle Ax,x\Bigr\rangle-\Bigl\langle x,u\Bigr\rangle\right)
\end{equation}
where $u\in H$ is  a fix point and $A$ is a mapping on $H$ which could be monotone, strongly
monotone or even bounded linear strongly positive operator as the case may be.

Let $K_1,K_2, \cdots ,K_N$ be $N$ closed convex subsets of a real Hilbert space $H$ having a
nonempty intersection $K.$ Suppose also that each $C_i$ is a fixed point set of nonexpansive
mappings $T_i : H\to H,$ $i = \overline{1,N}.$ Xu \cite{[47]} proved strong convergence of
the iterative algorithm
\begin{equation}
x_{n+1}=\alpha_nu + (I-\alpha_nA)T_nx_n, \quad  n\in \bn,
\end{equation}
where, $I$ is an identity mapping and $T_n := T_{n(mod N)},$ to a
unique solution of the quadratic minimization problem
\[\min\limits_{x\in K}\left(\frac12\Bigl\langle Ax,x\Bigr\rangle-\Bigl\langle x,u\Bigr\rangle\right),\]
where $A$ is a bounded linear strongly positive operator on $H$ and $u$ is a given point
in $H.$

Marino and Xu \cite{[26]}, proved that the iteration scheme given by
\begin{equation}\label{xnalphangammaATxn}
x_{n+1}=\alpha_n\gamma f(x_n) + (I-\alpha_nA)Tx_n, \quad  n\in \bn,
\end{equation}
converges strongly to a unique solution $x{'}\in F(T)$ of the variational inequality
\begin{equation}
\Bigl\langle(\gamma f-A)x{'},y-x{'}\Bigr\rangle\le 0, \quad \forall y\in F(T),
\end{equation}
which is the optimality condition for the minimization problem
\[\min\limits_{x\in F(T)}\left(\frac12\Bigl\langle Ax,x\Bigr\rangle-h(x)\right),\]
where $h$ is a potential function for $\gamma f$ (that is, $h'(x)=\gamma f(x)$ for all $x\in H$);
provided $f : H\to H$ is a contraction, $T : H\to H$ is nonexpansive and the iterative
parameter $\{\alpha_n\}$ satisfies appropriate conditions.

In \cite{[49]}, Yamada introduced the following hybrid iterative method
\begin{equation}\label{xnIalphanATxn}
x_{n+1}=(I-\alpha_n\mu A)Tx_n, \quad n\in \bn,
\end{equation}
where $T$ is nonexpansive, $A$ is $L-$Lipschitzian and strongly monotone operator with
constant $\eta$ and $0<\mu<\frac{2\eta}{L^2}.$ He proved that if $\{\alpha_n\}$ satisfies appropriate conditions,
then $\{x_n\}_{n=1}^\infty$ converges strongly to a unique solution $x{'}\in F(T)$ of the variational
inequality
\[\Bigl\langle Ax{'},y-x{'}\Bigr\rangle\ge 0, \quad \forall y\in F(T).\]

Recently, motivated by the iteration schemes \eqref{xnalphangammaATxn} and \eqref{xnIalphanATxn},
M. Tian \cite{[42]} introduced the following iterative method:
\begin{equation}
x_{n+1}=\alpha_n\gamma f(x_n)+(I-\alpha_n\mu A)Tx_n, \quad n\in \bn.
\end{equation}
Tian \cite{[42]} proved that if $f : H \to H$ is a contraction, $A : H\to H$ is an $\eta-$strongly
monotone mapping, $T : H\to H$ a nonexpansive mapping and the parameter  $\{\alpha_n\}$
satisfies appropriate conditions, then the sequence $\{x_n\}_{n=1}^\infty$ converges
strongly to a unique solution $x{'}\in F(T)$ of the variational inequality
\[\Bigl\langle(\gamma f-\mu A)x{'},y-x{'}\Bigr\rangle\le 0, \quad \forall y\in F(T).\]
His results generalized and improved the corresponding results of Marino and Xu
\cite{[26]} and Yamada \cite{[49]}.

In this paper, we extend  Tian's results \cite{[42]} to a finite family of nonexpansive mappings.
More precisely, we consider the following iterative algorithm
\begin{equation}\label{xn1alphangammafmuATxn}
x_{n+1}=\alpha_n\gamma f(x_n)+(I-\alpha_n\mu A)T_nx_n, \quad n\in \bn,
\end{equation}
where $I:H\to H$ is an identity mapping, $A:H\to H$ is a $k-$Lipschitzian $\eta-$strongly monotone mapping, $\{T_i\}_{i=1}^N:H\to H$ are nonexpansive mappings with $F:=\bigcap\limits_{i=1}^{N}F(T_i)\neq\emptyset,$ and $T_n:=T_{n(mod N)}.$ Under appropriate conditions, we shall prove strongly convergence of the sequence $\{x_n\}_{n=1}^\infty$ to a unique solution of the variational inequality
\[\Bigl\langle(\gamma f-\mu A)x{'},y-x{'}\Bigr\rangle\le 0, \quad \forall y\in F.\]
Our results generalize the corresponding results of Marino
and Xu \cite{[26]}, Tian \cite{[42]}, Xu \cite{[47]}, Yamada \cite{[49]}.

\section{Preliminaries}

The following lemmas play an important role in proving our main
results.

\begin{lemma}\cite{[47]}\label{Xusequence} Let $a_n$ be a sequence of nonnegative numbers satisfying
the following condition
\[a_{n+1}\le (1-\alpha_n)a_n+\alpha_n\beta_n\]
where $\{\alpha_n\}_{n=1}^\infty,$ $\{\beta_n\}_{n=1}^\infty$  are sequences of real numbers such that:
\begin{itemize}
  \item [(i)] $0<\alpha_n < 1$ and $\lim\limits_{n\to\infty}\alpha_n=0;$
  \item [(ii)] $\sum\limits_{n=1}^{\infty}\alpha_n=\infty;$
  \item [(iii)] $\limsup\limits_{n\to\infty}\beta_{n}\le 0.$
\end{itemize}
Then $\lim\limits_{n\to\infty}a_n=0$
\end{lemma}
\begin{lemma}\cite{[47]}\label{Xuvector}
In a real Hilbert space $H,$ there holds the following inequality
\[\|x+y\|^2\le \|x\|^2+2\langle y,x+y\rangle,\]
for all $x,y\in H.$
\end{lemma}
\begin{lemma}\cite{GK2}\label{GeobelKirknonexpansive} Let $H$ be a Hilbert space and $T: H\to H$ be a nonexpansive mapping with
$Fix(T)\neq\emptyset.$ If $\{x_n\}_{n=1}^\infty$ converges weakly to $x$ and $\{\|x_n-Tx_n\|\}_{n=1}^\infty$ converges to $0$, then $Tx=x.$
\end{lemma}
\begin{lemma}\cite{[42]}\label{Tianvarinequality}
Let $H$ be a real Hilbert space. Let $I : H \to H$ be an identity mapping, $f : H \to H$ be a contraction mapping with $0<\alpha<1$, $A:H\to H$ be a $k-$Lipschitzian $\eta-$strongly monotone mapping,  and $T:H\to H$ be a nonexpansive mapping with $F(T)\neq\emptyset.$ Let $x_t$ be a unique solution of the following equation
\[x_t=t\gamma f(x_t)+(I-t\mu A)Tx_t,\]
where $0<\mu<\frac{2\eta}{k^2},$ $\tau:=\mu(\eta-\frac{\mu k^2}{2}),$ $0<\gamma<\frac{\tau}{\alpha},$ and $0<t<1.$ Then the following variational inequality
\begin{equation}\label{variatinequality}
\Bigl\langle(\gamma f-\mu A)x,y-x\Bigr\rangle\le 0, \quad \forall y\in F(T).
\end{equation}
has a unique solution $x{'}$ on the set $F(T)$ and $x_t$ converges strongly to $x^{'}$ as $t\to 0.$
\end{lemma}

\section{Main results}

In this section we shall prove our main results. To formulate ones,
we need some auxiliary results.

\begin{lemma}
Let $H$ be a real Hilbert space. Let $I : H \to H$ be an identity mapping, $f : H \to H$ be a contraction mapping with $0<\alpha<1$, $A:H\to H$ be a $k-$Lipschitzian $\eta-$strongly monotone mapping,  and $\{T_i\}_{i=1}^N:H\to H$ be nonexpansive mappings with $F(T_i)\neq\emptyset,$ for all $i=\overline{1,N}.$ If $0<\mu<\frac{2\eta}{k^2}$  and $0<\alpha_n<1$ for all $n\in \bn$ then we have the following inquality
\begin{eqnarray}
\|(I-\alpha_n\mu A)T_nx-(I-\alpha_n\mu A)T_ny\|\le (1-\alpha_n\tau)\|x-y\|, \quad \forall n\in\bn,
\end{eqnarray}
for all $x,y\in H,$ where, $T_n:=T_{n(mod N)}$ and $\tau:=\mu(\eta-\frac{\mu k^2}{2}).$
\end{lemma}
\begin{pf} Let us denote by $S_n:=(I-\alpha_n\mu A)T_n.$ Then, we have for all $x,y\in H$ and $n\in\bn$
\begin{eqnarray*}
\|S_nx-S_ny\|^2&\le& \|T_nx-T_ny\|^2-2\alpha_n\mu\bigl\langle T_nx-T_ny,AT_nx-AT_ny\Bigr\rangle\\
&&+\alpha_n^2\mu^2\|AT_nx-AT_ny\|^2\\
&\le&\left(1-2\alpha_n\mu\eta+\alpha_n^2\mu^2k^2\right)\|T_nx-T_ny\|^2\\
&\le&\left(1-2\alpha_n\mu\left(\eta-\frac{\mu k^2}{2}\right)\right)\|x-y\|^2\\
&=&(1-2\alpha_n\tau)\|x-y\|^2\\
&\le&(1-\alpha_n\tau)^2\|x-y\|^2,
\end{eqnarray*}
which means
\begin{eqnarray}\label{SnxSnyconstnxy}
\|S_nx-S_ny\|\le (1-\alpha_n\tau)\|x-y\|.
\end{eqnarray} This completes the proof.
\end{pf}

\begin{theorem}
Let $H$ be a real Hilbert space. Let $I : H \to H$ be an identity mapping, $f : H \to H$ be a contraction mapping with $0<\alpha<1$, $A:H\to H$ be a $k-$Lipschitzian $\eta-$strongly monotone mapping,  and $\{T_i\}_{i=1}^N:H\to H$ be nonexpansive mappings with $F:=\bigcap\limits_{i=1}^NF(T_i)\neq\emptyset.$ Suppose that $0<\mu<\frac{2\eta}{k^2},$ $\tau:=\mu(\eta-\frac{\mu k^2}{2}),$ $0<\gamma<\frac{\tau}{\alpha},$ and the sequence $\{\alpha_n\}_{n=1}^\infty\subset (0,1)$ satisfies the following condition:
\begin{itemize}
  \item [(i)] $\lim\limits_{n\to\infty}\alpha_n=0;$
  \item [(ii)] $\sum\limits_{n=1}^\infty\alpha_n=\infty;$
  \item [(iii)] $\lim\limits_{n\to\infty}\dfrac{\alpha_n}{\alpha_{n+N}}=1.$
\end{itemize}
If $F=F(T_NT_{N-1}\cdots T_1)=F(T_1T_{N}T_{N-1}\cdots T_2)=F(T_2T_{1}T_NT_{N-1}\cdots T_3)=\cdots=F(T_{N-1}T_{N-2}\cdots T_1T_N)$ then the sequence $\{x_n\}_{n=1}^\infty$ which is defined by \eqref{xn1alphangammafmuATxn}, converges strongly to a unique solution $x{'}\in F$ of the variational inequality
\begin{eqnarray}\label{varinequalitygfmuA}
\Bigl\langle(\gamma f-\mu A)x{'},y-x{'}\Bigr\rangle\le 0, \quad \forall y\in F.
\end{eqnarray}
\end{theorem}
\begin{pf}
Since the mapping $T_NT_{N-1}\cdots T_1:H\to H$ is nonexpansive, then due to Lemma \ref{Tianvarinequality}, the variational inequality \eqref{varinequalitygfmuA} has a unique solution $x{'}$ on the set $F.$ We will show that the sequence $\{x_n\}_{n=1}^\infty$ given by \eqref{xn1alphangammafmuATxn} converges strongly to $x{'}.$

{\sc Step 1}. The sequences $\{x_n\}_{n=1}^\infty,$ $\{f(x_n)\}_{n=1}^\infty,$ $\{T_nx_n\}_{n=1}^\infty,$ and $\{AT_nx_n\}_{n=1}^\infty,$ are bounded. Indeed, let $p\in F$ and $S_n:=(I-\alpha_n\mu A)T_n.$ Using \eqref{SnxSnyconstnxy}, we then have
\begin{eqnarray}\label{xn1pxnp}
\|x_{n+1}-p\|&=&\|\alpha_n\gamma f(x_n)+S_nx_n-p\|\nonumber\\
&=&\|\alpha_n\gamma f(x_n)-\alpha_n\mu Ap+S_nx_n-S_np\|\nonumber\\
&\le&\alpha_n\gamma\alpha\|x_n-p\|+\alpha_n\|\gamma f(p)-\mu Ap\|+(1-\alpha_n\tau)\|x_n-p\|\nonumber\\
&=&(1-\alpha_n(\tau-\gamma\alpha))\|x_n-p\|+\alpha_n(\tau-\gamma\alpha)\frac{\|\gamma f(p)-\mu Ap\|}{\tau-\gamma\alpha}.
\end{eqnarray}
Since $0<\gamma<\frac{\tau}{\alpha},$ $0<\alpha_n<1,$ and $\lim\limits_{n\to\infty}\alpha_n=0$ then there exists $n_0>0$ such that $0<\alpha_n(\tau-\gamma\alpha)<1$ for all $n\ge n_0.$ We then obtain
\begin{eqnarray}
\|x_{n+1}-p\|\le \max\left\{\|x_{n}-p\|,\frac{\|\gamma f(p)-\mu Ap\|}{\tau-\gamma\alpha}\right\},\quad \forall n\ge n_0.
\end{eqnarray}
Therefore, we get
\begin{eqnarray}
\|x_{n+1}-p\|\le \max\left\{\|x_{n_0}-p\|,\frac{\|\gamma f(p)-\mu Ap\|}{\tau-\gamma\alpha}\right\},\quad \forall n\ge n_0,
\end{eqnarray}
which means, $\{x_n\}_{n=1}^\infty$ is a bounded sequence.

From the following inequality
\begin{eqnarray*}
&&\|f(x_n)-f(p)\|\le\alpha\|x_n-p\|,\\[2mm]
&&\|T_nx_n-p\|=\|T_nx_n-T_np\|\le\|x_n-p\|,\\[2mm]
&&\|AT_nx_n-Ap\|=\|AT_nx_n-AT_np\|\le k\|T_nx_n-T_np\|\le k\|x_n-p\|,
\end{eqnarray*}
it follows that the sequences $\{f(x_n)\}_{n=1}^\infty,$ $\{T_nx_n\}_{n=1}^\infty,$ $\{AT_nx_n\}_{n=1}^\infty,$ are bounded.

{\sc Step 2}. One has $\lim\limits_{n\to\infty}\|x_{n+1}-T_nx_n\|=0.$ Indeed, from $\lim\limits_{n\to\infty}\alpha_n=0,$ \eqref{xn1alphangammafmuATxn}, and Step 1 it follows that
\begin{eqnarray*}
\lim\limits_{n\to\infty}\|x_{n+1}-T_nx_n\|=\lim\limits_{n\to\infty}\alpha_n\|\gamma f(x_n)-\mu AT_nx_n\|=0
\end{eqnarray*}

{\sc Step 3}. For the sequence $\{x_n\}$ we have $\lim\limits_{n\to\infty}\|x_{n+N}-x_n\|=0.$ Indeed,  noting that $T_{n+N-1}=T_{n-1}$
 and $$S_{n+N-1}x_{n-1}-S_{n-1}x_{n-1}=\mu(\alpha_{n+N-1}-\alpha_{n-1})AT_{n-1}x_{n-1}$$ we then obtain
\begin{eqnarray*}
\|x_{n+N}-x_{n}\|&\le&\|\alpha_{n+N-1}\gamma f(x_{n+N-1})-\alpha_{n-1}\gamma f(x_{n-1})\|\\
&&+\|S_{n+N-1}x_{n+N-1}-S_{n-1}x_{n-1}\|\\
&\le& \alpha_{n+N-1}\gamma\alpha\|x_{n+N-1}-x_{n-1}\|+\gamma|\alpha_{n+N-1}-\alpha_{n-1}|\|f(x_{n-1})\|\\
&&+\|S_{n+N-1}x_{n+N-1}-S_{n+N-1}x_{n-1}\|\\
&&+\mu|\alpha_{n+N-1}-\alpha_{n-1}|\|AT_{n-1}x_{n-1}\|\\
&\le&(1-\alpha_{n+N-1}(\tau-\gamma\alpha))\|x_{n+N-1}-x_{n-1}\|\\
&&+|\alpha_{n+N-1}-\alpha_{n-1}|\Bigl(\gamma\|f(x_{n-1})\|+\mu\|AT_{n-1}x_{n-1}\|\Bigr)\\
&=&(1-\alpha_{n+N-1}(\tau-\gamma\alpha))\|x_{n+N-1}-x_{n-1}\|\\
&& +\alpha_{n+N-1}(\tau-\gamma\alpha)\beta_{n+N-1},
\end{eqnarray*}
where
\[\beta_{n-N+1}:=\frac{|\alpha_{n+N-1}-\alpha_{n-1}|}{\alpha_{n+N-1}}\cdot
\frac{\gamma\|f(x_{n-1})\|+\mu\|AT_{n-1}x_{n-1}\|}{\tau-\gamma\alpha}.\]
Letting $a_{n+N}:=\|x_{n+N}-x_{n}\|,$ one then gets
\begin{eqnarray}
a_{n+N}\le (1-\alpha_{n+N-1}(\tau-\gamma\alpha))a_{n+N-1}+ \alpha_{n+N-1}(\tau-\gamma\alpha)\beta_{n+N-1},
\end{eqnarray}
where
\begin{eqnarray*}
\limsup\limits_{n\to\infty}\beta_{n}=\limsup\limits_{n\to\infty}\frac{|\alpha_{n+N-1}-\alpha_{n-1}|}{\alpha_{n+N-1}}\cdot
\frac{\gamma\|f(x_{n-1})\|+\mu\|AT_{n-1}x_{n-1}\|}{\tau-\gamma\alpha}=0
\end{eqnarray*}
According to Lemma \ref{Xusequence}, we obtain
\[\lim\limits_{n\to\infty}a_{n+N}=\lim\limits_{n\to\infty}\|x_{n+N}-x_n\|=0.\]

{\sc Step 4.} We have $\lim\limits_{n\to\infty}\|x_n-T_{n+N-1}T_{n+N-2}\cdots T_{n}x_n\|=0.$ Indeed, noting that $T_n$ is a nonexpansive mapping and using Step 2 one has the following
\begin{eqnarray*}
\|x_{n+N}-T_{n+N-1}x_{n+N-1}\|&\to& 0,\\
\|T_{n+N-1}x_{n+N-1}-T_{n+N-1}T_{n+N-2}x_{n+N-2}\|&\to& 0,\\
\|T_{n+N-1}T_{n+N-2}x_{n+N-2}-T_{n+N-1}T_{n+N-2}T_{n+N-3}x_{n+N-3}\|&\to& 0,\\
&\vdots&\\
\|T_{n+N-1}T_{n+N-2}\cdots T_{n+1}x_{n+1}-T_{n+N-1}T_{n+N-2}T_{n+N-3}\cdots T_{n}x_{n}\|&\to& 0,
\end{eqnarray*}
as $n\to\infty.$ Hence, we obtain
\[\|x_{n+N}-T_{n+N-1}T_{n+N-2}T_{n+N-3}\cdots T_{n}x_{n}\|\to 0,\quad n\to\infty.\]
Since $\|x_{n+N}-x_n\|\to 0$ as $n\to\infty$ (see Step 3) we then get
\[\lim\limits_{n\to\infty}\|x_n-T_{n+N-1}T_{n+N-2}\cdots T_{n}x_n\|=0.\]

{\sc Step 5}. We want to show that $\limsup\limits_{n\to\infty}\Bigl\langle(\gamma f-\mu A)x{'}, x_{n}-x{'}\Bigr\rangle\le 0,$ where $x{'}$ is a unique solution of the variational inequality \eqref{varinequalitygfmuA} for the nonexpansive mapping $T_NT_{N-1}\cdots T_1:H\to H$ on the set $F.$

Let $\{x_{n_k}\}_{k=1}^\infty$ be a subsequence of $\{x_n\}_{n=1}^\infty$ such that
\[\limsup\limits_{n\to\infty}\Bigl\langle(\gamma f-\mu A)x{'}, x_{n}-x{'}\Bigr\rangle=\lim\limits_{k\to\infty}\Bigl\langle(\gamma f-\mu A)x{'}, x_{n_k}-x{'}\Bigr\rangle.\]
Since $\{x_{n_k}\}_{k=1}^\infty$ is  bounded and $H$ is a real Hilbert space, then there exists a subsequence $\{x_{n_{k_m}}\}_{m=1}^\infty$ of $\{x_{n_k}\}_{k=1}^\infty$ such that $\{x_{n_{k_m}}\}_{m=1}^\infty$ converges weakly to some point $y\in H.$ Without loss any generality, we may assume that $n_{k_m}$ are such kind of numbers that $T_{n_{k_m}}=T_{i_0}$ for some $i_0\in\{1,2,\cdots, N\}$ and for all $m\in \bn.$ Then, from Step 4, it follows that
\begin{eqnarray*}
\lim\limits_{n\to\infty}\|x_{n_{k_m}}-T_{i_0+N-1}T_{i_0+N-2}\cdots T_{i_0}x_{n_{k_m}}\|=0.
\end{eqnarray*}
So due to Lemma \ref{GeobelKirknonexpansive}, we have $y\in F(T_{i_0+N-1}T_{i_0+N-2}\cdots T_{i_0})=F.$
Therefore,
\begin{eqnarray*}
\limsup\limits_{n\to\infty}\Bigl\langle(\gamma f-\mu A)x{'}, x_{n}-x{'}\Bigr\rangle&=&\lim\limits_{m\to\infty}\Bigl\langle(\gamma f-\mu A)x{'}, x_{n_{k_m}}-x{'}\Bigr\rangle\\
&=&\Bigl\langle(\gamma f-\mu A)x{'}, y-x{'}\Bigr\rangle\le 0.
\end{eqnarray*}

{\sc Step 6}. The sequence $\{x_n\}$ converges to $x'$, i.e. $\lim\limits_{n\to\infty}\|x_{n}-x{'}\|=0.$  Indeed, using Lemma \ref{Xuvector} and \eqref{xn1pxnp} we obtain
\begin{eqnarray*}
\|x_{n+1}-x{'}\|^2&=&\|S_nx_n-S_nx{'}+\alpha_n(\gamma f(x_n)-\mu AT x{'})\|^2\\
&\le&\|S_nx_n-S_nx{'}\|^2+2\alpha_n\Bigl\langle\gamma f(x_n)-\mu Ax{'}, x_{n+1}-x{'}\Bigr\rangle\\
&\le&(1-\alpha_n\tau)^2\|x_n-x^{'}\|^2+2\alpha_n\gamma\Bigl\langle f(x_n)-f(x{'}), x_{n+1}-x{'}\Bigr\rangle\\
&&+2\alpha_n\Bigl\langle(\gamma f(x_n)-\mu A) x{'}, x_{n+1}-x{'}\Bigr\rangle\\
&\le&(1-\alpha_n\tau)^2\|x_n-x^{'}\|^2+2\alpha_n\gamma\alpha\|x_n-x{'}\|\|x_{n+1}-x{'}\|\\
&&+2\alpha_n\Bigl\langle(\gamma f(x_n)-\mu A) x{'}, x_{n+1}-x{'}\Bigr\rangle\\
&\le&(1-\alpha_n\tau)^2\|x_n-x^{'}\|^2+2\alpha_n(1-\alpha_n(\tau-\gamma\alpha))\gamma\alpha\|x_n-x{'}\|^2\\
&&+2\alpha_n^2\gamma\alpha\|\gamma f(x{'})-\mu Ax{'}\|+2\alpha_n\Bigl\langle(\gamma f(x_n)-\mu A) x{'}, x_{n+1}-x{'}\Bigr\rangle\\
&\le&(1-2\alpha_n(\tau-\gamma\alpha))\|x_n-x^{'}\|^2+2\alpha_n(\tau-\gamma\alpha)\beta_n,\\
\end{eqnarray*}
where
\begin{eqnarray*}
\beta_n&=&\frac{\Bigl\langle(\gamma f(x_n)-\mu A) x{'}, x_{n+1}-x{'}\Bigr\rangle}{\tau-\gamma\alpha}\\
&&+\alpha_n\left[\left(\frac{\tau^2}{2(\tau-\gamma\alpha)}-\gamma\alpha\right)\|x_n-x{'}\|^2+\frac{\gamma\alpha}{\tau-\gamma\alpha}\|\gamma f(x{'})-\mu Ax{'}\|\right].
\end{eqnarray*}
Letting $a_n:=\|x_n-x{'}\|^2,$ we then have
\begin{eqnarray}
a_{n+1}\le (1-2\alpha_n(\tau-\gamma\alpha))a_n+2\alpha_n(\tau-\gamma\alpha)\beta_n.
\end{eqnarray}
Since $\limsup\limits_{n\to\infty}\Bigl\langle(\gamma f-\mu A)x{'}, x_{n+1}-x{'}\Bigr\rangle\le 0$ (see Step 5) one can get
\begin{eqnarray*}
\limsup\limits_{n\to\infty}\beta_n\le 0.
\end{eqnarray*}
Therefore, according to Lemma \ref{Xusequence}, we obtain that
\[\lim\limits_{n\to\infty}a_n=\lim\limits_{n\to\infty}\|x_n-x{'}\|^2=0,\]
which implies that $\lim\limits_{n\to\infty}\|x_n-x{'}\|=0.$ This completes the proof.
\end{pf}

Note that if $f(x)=\gamma u$ for some $\gamma\in (0,1)$ and $u\in H$ we recover Xu's \cite{[47]} result.  Our results also generalize the corresponding results of Marino
and Xu \cite{[26]}, Tian \cite{[42]}, Yamada \cite{[49]}.

\section*{Acknowledgement}  This
work was done while the first named author (F.M.) was visiting the Abdus Salam International Centre for Theoretical Physics, Trieste,
Italy as a Junior Associate. He would like to thank the Centre for hospitality and financial support.

\end{document}